\documentclass[14pt,a4paper]{article}
\usepackage{Packages}

\begin{document}


\vspace{0.5 cm}

{\centering

{\bfseries\large Opportunities for Hybrid Modeling Approaches in Energy Systems optimization} 


\vspace{0.5 cm}
Mohamed Tahar Mabrouk$^{1*}$, Shri Balaji Padmanabhan$^{1}$, Bruno Lacarrière$^{1}$, Benoit Delinchant$^{2}$, Sacha Hodencq$^{2}$, Xavier Roboam$^{3}$, Bruno Sareni$^{3}$, Mathieu Vallee$^{4}$\\ [0.2 cm]

$^{1}$IMT Atlantique, Department of Energy Systems and Environment, GEPEA UMR CNRS, 6144, F-44307, Nantes, France\\ [0.2 cm]
$^{2}$University of Grenoble Alpes, Grenoble INP, Grenoble, France\\ [0.2 cm]
$^{3}$ENSEEIHT, Université de Toulouse, Toulouse, France\\ [0.2 cm]
$^{4}$Univ. Grenoble Alpes, CEA, Liten, Campus Ines, 73375 Le Bourget du Lac, France\\ [0.2 cm]
\textit{*Corresponding Author: mohamed-tahar.mabrouk@imt-atlantique.fr}\\ 
\par}
\vspace{0.5 cm}

{\bfseries \centering{ABSTRACT} \par }

This paper surveys the primary computational hurdles of Energy Systems optimization coming from different sources: model-induced complexity, optimization algorithm requirements, and uncertainties handling (both aleatoric and epistemic). Techniques to reduce complexity such as time-series and spatial aggregation, model order reduction, and specialized optimization strategies are reviewed for their effectiveness in balancing computational feasibility and model fidelity. Furthermore, Various uncertainty-management frameworks, including scenario-based approaches, robust optimization, and distributionally robust methods, are reviewed and their limitations in scaling and data requirements are discussed. The potential of hybrid modeling emerges as a key avenue: by fusing mechanistic and machine learning elements, hybrid techniques for modelling and optimization can harness the strengths of both worlds while mitigating their respective drawbacks. The paper highlights several directions for further research to develop advanced methods to tackle the complexity of MES.

\section{INTRODUCTION}
In light of the ongoing energy transition aiming for higher efficiency, reduced environmental impact (\cite{noauthor_decisive_2021}), and the growing integration of intermittent renewable energy sources (\cite{i2024a}), designing highly efficient and reliable energy systems that exploit the synergy between energy sectors and flexibility sources while maintaining acceptable Quality of Service (QoS) has become crucial. Future energy systems will incorporate storage devices, energy conversion technologies, and will implement complex demand response mechanisms (\cite{i2024a}). Achieving optimal design and operation of such systems requires accurate modeling of system behaviors, combined with efficient optimization algorithms that identify the best decisions.

Energy research generally falls into two distinct domains (\cite{subramanian_modeling_2018}): In one hand, there is the technological domain. In the other hand there is the Energy Policy (EP) domain which includes economy and social science. Each of them is characterized by different objectives, modeling strategies, and levels of technological detail. Technological models concentrate on unit operations, individual plants, or supply chains, treating technological parameters endogenously and incorporating economic, environmental, or social factors as exogenous inputs when needed. These models typically support design, operations, and control decisions. Conversely, EP models aggregate the entire energy sector at regional, national, or global scales, embedding economic variables (e.g., supply, demand, and market equilibrium) within the model while incorporating technology and other factors as exogenous components. Traditionally used for high-level, long-term planning, EP models now increasingly address operational decision-making by incorporating short-term elements such as seasonal demand, price spikes, and weather forecasts. This is made necessary by the expansion of intermittent renewable energy sources and the progressive deregulation of electricity and energy markets.
Energy system modeling must carefully adapt its level of detail to both the intended application (control, operation, design, and planning) and the size of the system, which can range from a small facility to large, cross-border networks. Highly detailed models may provide fine predictions on the system's and technologies performance but can become computationally intractable at larger scales and highly uncertain in pre-design phases, whereas more aggregated models may overlook essential system behaviors preventing informed decision-making. Striking a balance between accuracy and computational feasibility is therefore essential to maintain a tractable, model-based decision-making process.

Optimization algorithms are used along with the system's model to evaluate different solutions against objectives like cost, emissions, or reliability, systematically searching for the best decision. However, as the system’s complexity and number of decision variables grow, the computational effort required by this approach can become prohibitively large, especially for highly detailed or large-scale energy systems. Moreover, incorporating uncertainties into the optimization process implies adding constraints to capture uncertainty or repeatedly solving variations of the optimization problem under different assumptions which further increases the complexity of the decision-making.

In this paper, we present an in-depth review of the existing literature on mechanistic, model-based methodologies for Energy Systems optimization. By systematically identifying and categorizing the challenges associated with these traditional approaches, we highlight recurring limitations in scalability, complexity handling, and uncertainty management across different decision-making levels and system sizes. Building on these insights, the paper proposes key directions for future research in hybrid modeling techniques, offering pathways to address problems that have remained intractable using traditional methods.

\section{Model-based energy systems optimization}
\subsection{Applications}
Model-based energy systems optimization aim to take decisions that minimize or maximize one or more Key Performance Indicators (KPIs) estimated using the system's model while complying with technical, economic, environmental and legal constraints imposed by stakeholders, as well as the physical limitations of the system and its environment.

\begin{figure}[h!]
    \centering
    \includegraphics[scale=0.7]{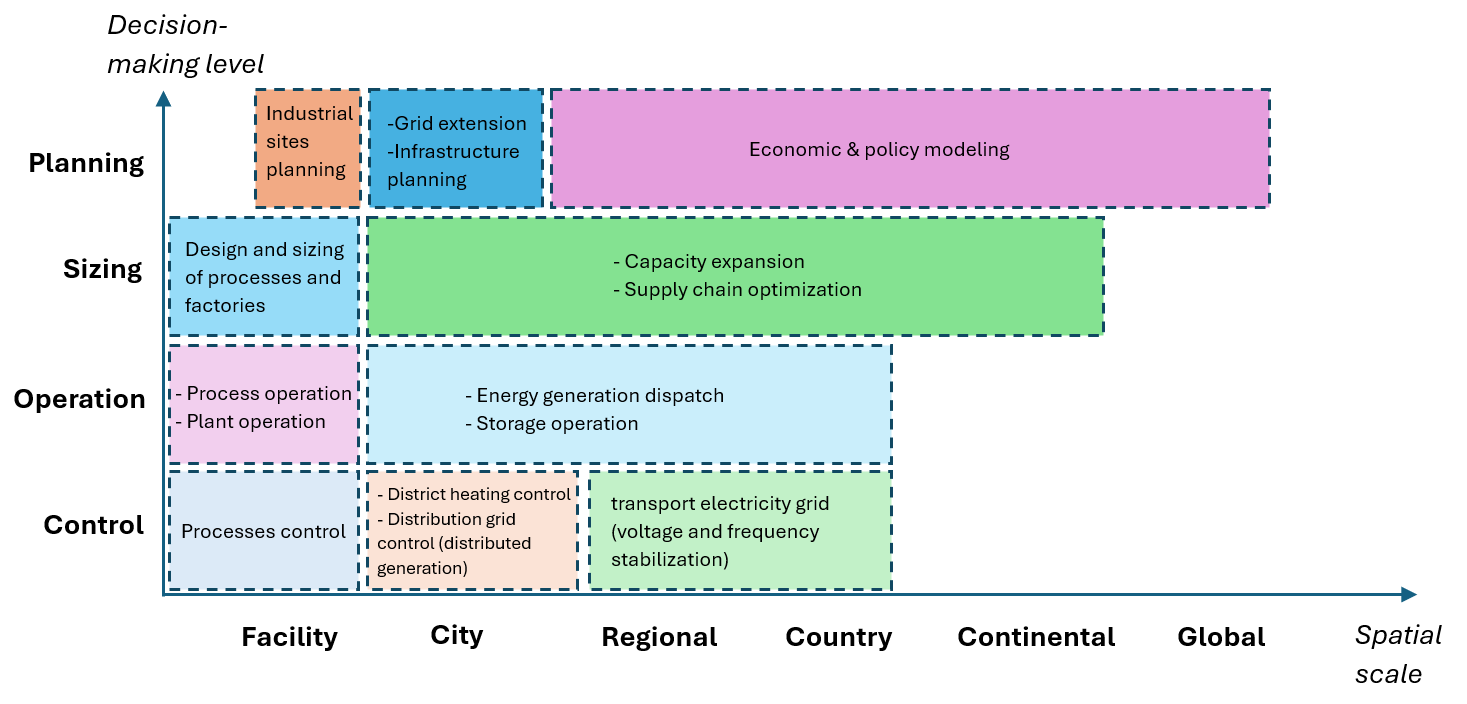}
    \caption{Decision making levels VS energy system scales, adapted from (\cite{WEIDNER2022107883})}
    \label{fig:decision levels}
\end{figure}

This optimization process applies to a wide range of applications, as shown in the figure \ref{fig:decision levels}. These applications are classified based on the decision-making level and system's scale. There are four main decision-making levels: \textbf{1. Control} which involves real-time decisions where actions are executed instantly by the system's actuators. Here, the system already exists, and some of its technical characteristics are known. \textbf{2. Operation} which focuses on planning short-term operation strategies, covering a time frame from a few hours to several months. \textbf{3. Sizing} which determines the size and technical specifications of the system's components to ensure efficient performance. \textbf{4. Planning} which deals with long-term decisions, such as selecting system components, determining their placement, and defining their interconnections or more broadly optimizing economic and political policy decisions. This framework applies to energy systems of various sizes, from individual components to the global scale. 
\subsection{Mathematical formulations and solving paradigms}
A decision-making process is a structured sequence of steps aiming to evaluate different possible decisions and selecting the option that best balances the criteria or maximizes the desired outcome. For simplification, in this paper will be abstracted using the general formulation of a mono-objective constrained optimization problem. Knowing that all the subsequent discussion could be naturally extended to any multi-objective decision-making process:
\begin{equation}
 \begin{aligned}
\min_{X \in \Gamma(X, \theta)} \quad  & f(X,\theta)\\
\end{aligned}
    \label{eq:optim_mathematical_formulation}
\end{equation}
where \( X \) represents the set of optimization variables, defined as \( X = X_d \cup X_s \), with \( X_d \) being the set of \textbf{decision variables}, and \( X_s \) being the set of \textbf{simulation variables}. The parameters set is denoted as \( \theta \), while \( f(X,\theta) \) represents the \textbf{objective function} to be minimized. The feasible region \( \Gamma(X, \theta) \) consists of a set of constraints, expressed as:

\begin{equation}
\Gamma(X, \theta) = \Gamma_b(X, \theta) \cup \Gamma_t(X, \theta) \cup \Gamma_s(X, \theta)
\end{equation}
where, \( \Gamma_b(X, \theta) \) are \textbf{boundary constraints} which encompass environmental, legal, and economic requirement, as well as stakeholder preferences, \( \Gamma_t(X, \theta) \) includes \textbf{technical constraints} which represent the system's own technical limitations, and \( \Gamma_s(X, \theta) \) corresponds to \textbf{simulation constraints}, i.e., the system model equations.

\begin{figure}[h!]
    \centering
    \includegraphics[scale=0.3]{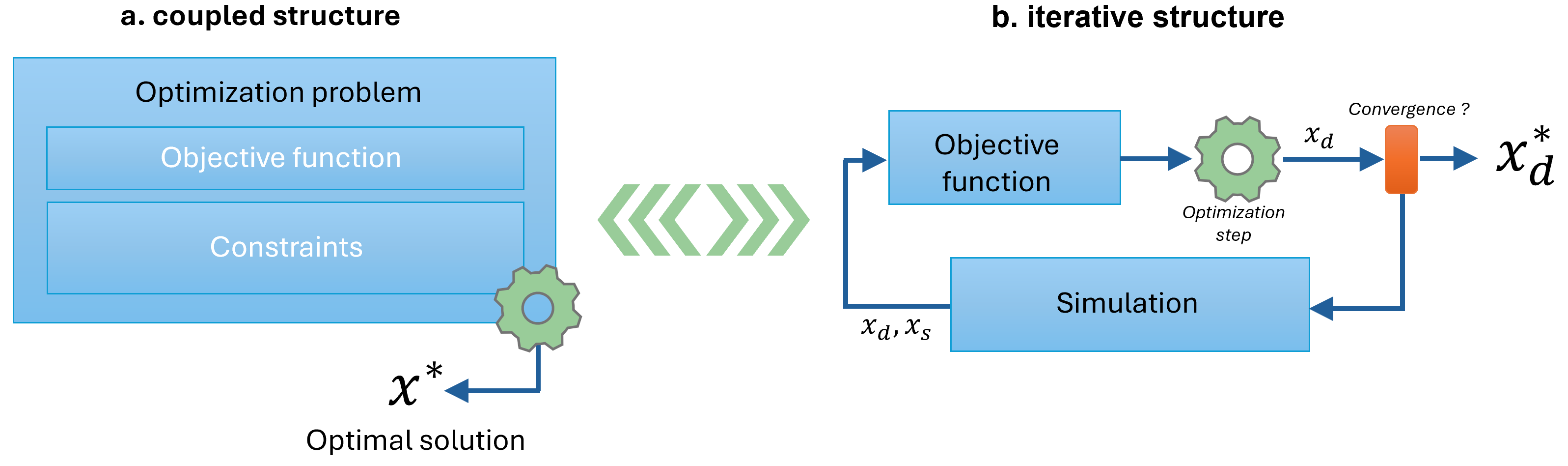}
    \caption{Optimization solving paradigms: coupled structure Vs. iterative structure}
    \label{fig:diagram}
\end{figure}

Figure \ref{fig:diagram} illustrates two distinct approaches to solving optimization problems: the \textbf{coupled structure} and the \textbf{iterative structure}. In the \textit{coupled structure}, the optimization problem is formulated as a single problem, where the objective function and constraints are embedded within the same framework. The solution process directly determines the optimal solution \( x^* \) in the same time, making this approach well-suited for problems with well-defined mathematical formulations that can be efficiently solved using optimization solvers. In contrast, the \textit{iterative structure} (right) decomposes the optimization process into an iterative loop consisting of an optimization step and a simulation step. Here, decision variables \( x_d \) are updated at each optimization step and simulation variables \( x_s \) are updated during the simulation step. The loop converges with to an optimal decision \( x^*_d \). This approach is particularly useful for complex systems. 

\section{Challenges}
Solving optimization problems for energy systems presents significant challenges that can be broadly categorized into two categories: computational complexity and uncertainty management.
\subsection{Algorithmic and computational complexity}
\label{sec: complexity}
\begin{figure}[h!]
    \centering
    \includegraphics[scale=0.35]{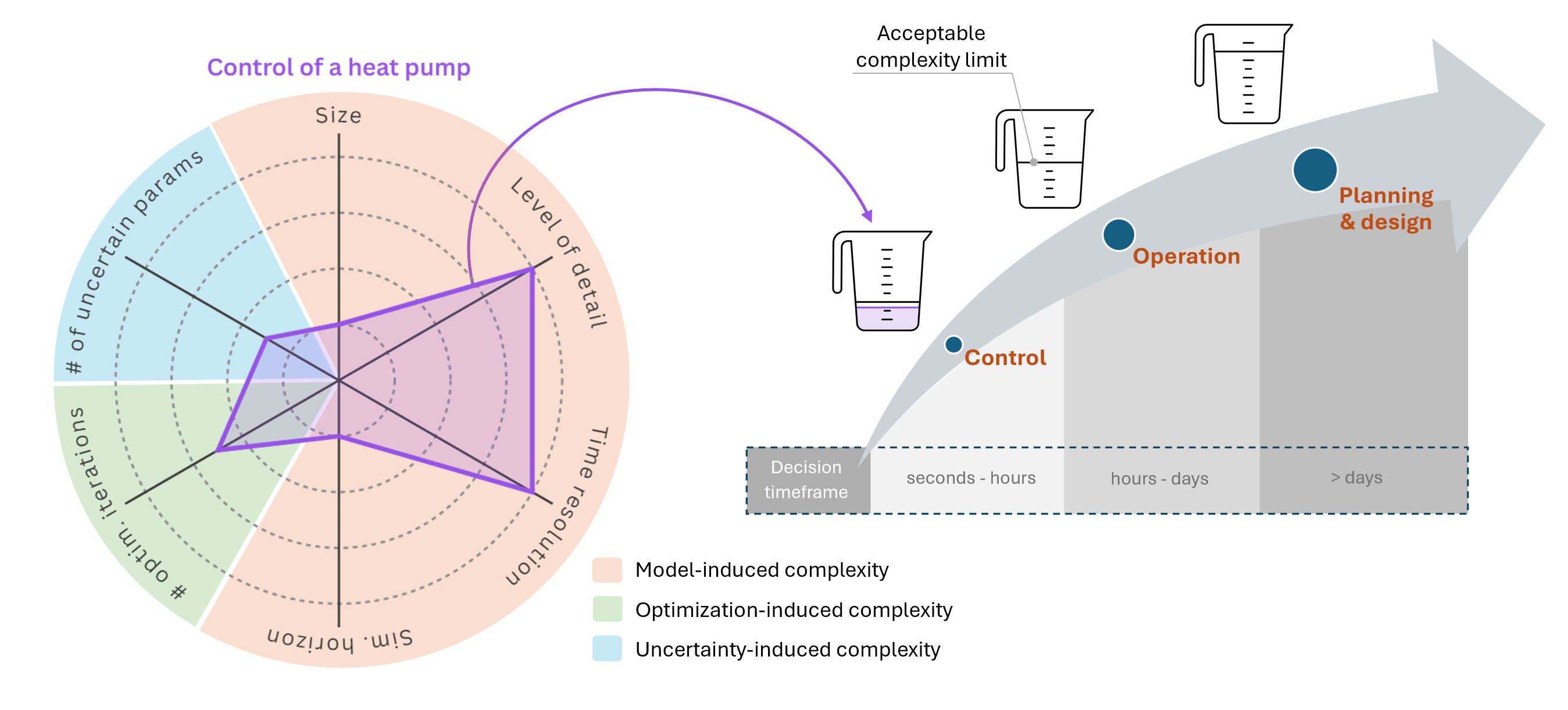}
    \caption{Model's tractability based on complexity determinants and the targeted application : example of the control of a heat pump}
    \label{fig:complexity determinants}
\end{figure}

As depicted in Figure (\ref{fig:complexity determinants}), different decision-making levels spans various time-frames, ranging from short-term control (seconds to hours), through medium-term operational decisions (hours to days), to long-term strategic planning and design decisions (days or longer). Each application inherently has its own acceptable complexity limit. The manageability of a decision-making process relies significantly on keeping its complexity below the complexity limit of the targeted application. Complexity in decision-making processes arises from three main categories: model-induced complexity, optimization-induced complexity, and uncertainty-induced complexity.

\subsubsection{Model-induced complexity} 
 The majority of optimization variables and constraints originate from the system's model. Larger systems and/or more detailed models increase the number of variables and constraints which substantially impacts the computational complexity of the optimization problem. 

By definition, "a model is a simplified representation of reality designed for a specific purpose" (\cite{kotzur_modelers_2021}), in a given context. The challenge resides in accurately capturing the essential aspects of reality while maintaining the lowest possible level of detail. A model must be detailed enough to provide useful insights but simplified enough to remain computationally manageable and applicable to its intended purpose. For instance, real-time control of a heat pump requires a high level of detail and high time resolution to accurately capture fast system's dynamics. However, due to strict computational constraints imposed by the short decision time-frame, the decision-making process must be simplified by adopting shorter simulation horizons and minimizing the level of detail as much as possible. Conversely, operational and strategic planning decisions typically allow greater complexity, enabling models to incorporate larger system scales and extended simulation horizons. Nonetheless, even for these applications, the level of detail (\cite{priesmann_are_2019}) and temporal resolution (\cite{poncelet_importance_2014}) should be carefully managed to balance complexity and accuracy. Identifying the appropriate complexity level tailored to the targeted application ensures computational feasibility and effective optimization in energy systems. However, there is non systematic and comprehensive approach to achieve this task.
\subsubsection{Optimization-induced complexity}
The choice of the optimization algorithm depends on whether we are looking for the global solution with or without guarantee or only a local solution (\cite{kochenderfer_algorithms_2019}). The complexity of the optimization algorithm adopted for the resolution process is determined by the number of iterations required by the solver to reach an optimal solution, and consequently, by the number of times the objective function and the constraints set must be evaluated. This complexity is tightly linked to the size of the model, specifically the number of variables and constraints included. Additionally, the nature of the objective function and constraints significantly influences the complexity of the optimization algorithm. Optimization techniques can be classified into convex and non-convex optimization (\cite{lin_review_2012}). Convex optimization problems are divided into linear and nonlinear optimization, for both of which convergence to a global optimum is guaranteed by state-of-the-art solvers. However, even within this simpler class, complexity can drastically increase with the number of variables; for instance, interior point algorithms which are widely used for linear programming (LP) have polynomial complexities (\cite{Ayache_2023}). Non-convex optimization problems are notably more challenging due to complexities arising from nonlinearities or discrete variables in constraints and/or the objective function. These characteristics create difficult optimization landscapes, often characterized by multiple local optima and problematic convergence behaviors.

\subsubsection{Uncertainty-induced complexity}
The presence of uncertain parameters adds significant complexity to decision-making by increasing the number of potential scenarios. This complexity depends on the method used to handle these uncertainties. These methods necessitate either solving the problem multiple times, evaluating the model for multiple scenarios, or increasing the number of constraints in the problem. Sources of uncertainty and methods used for managing them in energy-systems optimization will be discussed in detail in the following sections.
\subsection{Uncertainties}
\subsubsection{Aleatoric uncertainties}
Aleatoric uncertainty refers to the inherent randomness or variability in a process or system that cannot be eliminated by collecting more data or improving models. It is the uncertainty that comes from natural fluctuations of the systems environment. Sources of aleatoric uncertainties are multiple.
\textbf{Weather and Ambient Conditions:}
  Parameters such as temperature, humidity, precipitation, solar irradiation and wind speed fluctuate and affect  supply (renewables), demand (heating/cooling), and systems' efficiencies.
\textbf{  Load Forecasts:}
  Demand in residential, commercial, and industrial settings often deviates from predictions due to changing occupancy, process schedules, and other unforeseen shifts in consumption.
\textbf{Human Behavior:}
  User preferences and behavior can cause significant deviations from the expected performance. This includes among others comfort levels, occupancy patterns, and responsiveness to pricing.
\textbf{Market, costs and Price Volatility:}
  This includes Uncertain energy prices, components, labor and material costs. This introduces risks in both investment and operation decisions.
\textbf{Equipment Reliability and Outages:}
  Forced outages and faster-than-expected component degradation can occur and are hard to predict.
\textbf{Measurement Uncertainties:}
  Sensor inaccuracies, communication delays, and imperfect state estimation can obscure the true operating conditions. This is particularly challenging for real-time control.
\textbf{Regulatory and Policy Changes:}
Unexpected changes in emissions regulations, market structures, or incentive mechanisms can drastically influence operational strategies and capital planning.

Uncertain parameters can be either continuous or discrete, and the choice of how to represent uncertainty and handle them is often dictated by practical constraints on data availability.

\subsubsection{Epistemic uncertainties}
Epistemic uncertainty, also known as model-form uncertainty, refers to the systematic deviation between the predictions of a model and the behavior of the modeled system in the real world. It can emerge when a mechanistic model is simplified or when the underlying phenomena are not fully or accurately represented or understood (\cite{EUGENE2023108430}). Another definition, widely adopted by the energy planning community states that epistemic uncertainty is the absence of knowledge about the potential value of the uncertain parameters (\cite{rodriguez-matas_improving_2024}). Although different, this definition is actually a special case of the first one. Indeed, the absence of knowledge regarding certain parameter values implies uncertainty about whether these parameters should be treated as fixed constants or modeled as variables using some unknown models. Unlike aleatoric uncertainties, epistemic uncertainties can be explained away if we provide enough data. However, in practice, distinguishing between aleatoric and epistemic uncertainties is inherently challenging. For instance, during model calibration for predictive control applications, identifying whether discrepancies between observed data and model predictions arise from suboptimal parameter estimation or fundamental limitations in the model's structure requires a rigorous statistical framework.

\section{Classic approaches to tackle these challenges}
\label{sec:classic}
\subsection{Complexity reduction} 
\label{sec:complexity_reduction}

 
\subsubsection{Temporal reduction}
As demonstrated in section \ref{sec: complexity}, two temporal factors affect the complexity of the model: time resolution and the simulation horizon. A high temporal resolution is crucial in small time frame decision such as control or operation. However, even for longer term decision such as design and planning, a sufficient time resolution is important to capture renewable energy sources intermittency. Simulation horizon is more increasingly impactful for sizing and planning. The idea behind temporal reduction is to reduce either the temporal resolution or the simulation period of the model which leads to a reduction of the decision-making problem complexity. The use of one or both of these srategies can be referred in general to Time series Reduction (TSR).

Temporal resolution reduction can be achieved either by down-sampling or segmentation (\cite{kotzur_modelers_2021}). Down-sampling consists of regularly replacing a predefined number of adjacent time steps by one average time step. This tends to shaving peak values which leads to underestimating installed capacities and overestimating capacity factors (\cite{pfenninger_dealing_2017}). Segmentation is a more advanced technique consisting of merging time steps based on a metric of their similarity which lead to irregular new time steps. Many methods have been proposed in the literature including optimization (\cite{vom_stein_development_2017}), clustering (\cite{teichgraeber_clustering_2019}), or a combination of those (\cite{rigo-mariani_optimized_2022}). 

Decreasing the simulation period is achieved in the literature by grouping homogeneous periods (or time steps) together and replacing them by typical periods. A multitude of clustering techniques are used to perform this task including k-means, k-medoids, k-shape, DBA clustering, agglomerative hierarchical clustering, etc..(\cite{teichgraeber_clustering_2019}). In all the cases, the number of typical periods is predefined and highly affects the accuracy of the results. \cite{goke_adequacy_2022} compared loss of load and system cost of different TSR methods by applying the to a macro-energy system with high renewable sources penetration and extensive seasonal variations. They concluded that TSR must be applied with care, particularly because  seasonal storage requirements and intermittency can diminish their accuracy. They also highlighted the importance of refining TSR methods, for instance, by identifying extreme situations or using adaptive resolution. \cite{hoffmann_typical_2021} investigated the impact of different temporal aggregation techniques on the accuracy and computational complexity of two energy system models: a cogeneration system for a self-sufficient building and a multi-regional electricity dispatch model. The results underscore that choices around temporal aggregation depend more on storage cycle lengths than on autocorrelation patterns of input time series. Moreover, clustering indicators showed only a weak correlation with model accuracy once the aggregated time series were incorporated into the optimization. The authors conclude that understanding of a given energy system model’s mathematical structure is vital to selecting the most effective aggregation strategy.
\subsubsection{Spatial reduction}
Spatial reduction techniques aim to reducing the number of model's regions or components. These techniques are widely used in energy system planning studies where the scale of the system is typically large and potential decrease in computation times is substantial. 
\cite{bogdanov_hierarchical_2023} used Hierarchical optimization by solving the optimization problem with aggregated macro-regions (low resolution), solving the problem for each macro-region as an island and then reconstituting full results in the last step. This approach lead to $2.35\times$ to $3.3\times$ acceleration with a maximal deviation of $3\%$ compared  to the reference. 
 \cite{frysztacki_comparison_2022} compared different hierarchical clustering algorithms for the optimization of the European electrical system. All clustering algorithms tested exhibited comparable performances and any of them perform better than country boundaries-based clustering (\cite{brown_synergies_2018}) or k-means on the coordinates of the network (\cite{frysztacki_strong_2021}).
 \cite{patil_impact_2021} studied the impact of both regional and technology aggregation on the Total Annual Cost (TAC) of an instance of the European energy system. They found that 33 regions and 38 renewable technology types in each region is optimal to keep TAC deviation under $5\%$ and achieve a run time reduction by $93\%$. Finally, \cite{jacobson_quantifying_2024} used other metrics that TAC to compare optimization performance such as technology siting accuracy and CO2 emissions and showed that, while increasing resolution increases the prices accuracy asymptotically, the accuracy improvement for the other metrics is not asymptotic. They demonstrated also that spatial resolution is more impactful than temporal resolution because of the loss of information on the system's structure when lower spatial resolution is used.

\subsubsection{Model reduction}

\begin{table}[h]
\centering
\caption{Classification of ROM Methods by Intrusiveness, Data-Driven Approach, and Applicability to Linear/Nonlinear Systems}
\label{table: ROM_classification}
\begin{tabular}{p{0.4\textwidth}p{0.15\textwidth}p{0.15\textwidth}p{0.2\textwidth}}
\toprule
\textbf{Method} & \textbf{Intrusive \newline Non-Intrusive} & \textbf{Data-Driven?} & \textbf{Linear \newline Nonlinear?} \\
\midrule

\textbf{POD (Proper Orthogonal Decomposition)} & & & \\
\quad \textit{POD + Galerkin Projection} & Intrusive  & Yes & Linear or nonlinear \\
\quad \textit{POD + Interpolation} & Non-intrusive & Yes  & Linear or nonlinear \\
\addlinespace

\textbf{Balanced Truncation} & & & \\
\quad \textit{Classical Balanced Truncation} & Intrusive & No & Linear  \\
\quad \textit{Empirical Balanced Truncation} & Non-intrusive & Yes & Primarily linear; \newline nonlinear extensions \\
\addlinespace

\textbf{Krylov Subspace } & Intrusive & No & Typically linear \\
\addlinespace

\textbf{DMD (Dynamic Mode Decomposition)} & Non-intrusive  & Yes  & Primary linear \newline nonlinear extensions \\
\addlinespace

\textbf{System Identification} &  &  &  \\
\quad \textit{ERA (Eigensystem Realization Algorithm)} & Non-intrusive & Yes & Linear \\
\quad \textit{N4SID} & Non-intrusive & Yes & Linear \\
\quad \textit{Loewner Framework} & Non-intrusive & Yes & Linear \\
\addlinespace

\textbf{Neural Network Surrogates} & Non-intrusive  & Yes  & Linear or nonlinear\\
\addlinespace

\bottomrule
\end{tabular}
\end{table}

Model reduction is not relevant for large optimal planning of large systems where it is demonstrated that temporal and spatial aggregation are generally more effective at reducing complexity than simplifications in model structure \cite{priesmann_are_2019}. However, for applications such as optimal control, where capturing fine dynamic behavior is essential, or optimal design of technologies or components, where high-fidelity predictions are required, model reduction offers greater opportunities to simplify complexity while retaining essential details. 
Model reduction is used to find a low-order model that approximates the behavior of the original high-order model enabling fast computation while introducing a slight loss of accuracy (\cite{besselink_comparison_2013}).
The development of a Reduced Order Model (ROM) involves two main stages: Offline stage where the full Order Models is used with selected parameter values using HPC resources to extract features that will be used to build the surrogate model. Online stage where the surrogate model is used to predict solutions for new parameter values, offering significant computational speedup.
Methods used to develop ROMs are typically divided into intrusive methods, which involve manipulating the governing equations, and non-intrusive methods, which rely only on simulation data without using the original model's structure (\cite{padula_brief_2024}). All non-intrusive methods are by definition data-driven. However some intrusive methods are data-driven too for instance when simulation snapshots are needed to generate a reduced basis and the original model's equations are then projected in the new reduced space. Table \ref{table: ROM_classification} classifies some of the most used ROM techniques in the literature. 
In one hand, the main limitation of of intrusive ROMs are the necessity to manipulate directly the model's equations which limits their usage in real world applications where simulators are frequently complex and closed source. In the other hand, the limitations of non-intrusive data driven ROMs are multiple, including the necessity to run expensive simulations for a big number of parameters realizations. This is particularly challenging when the dimensionality of the parameter set is high, exploding the number of realizations needed to map all the possibilities (i.e. the curse of dimensionality) (\cite{mufti_multi-fidelity_2024}).

\subsubsection{Optimization strategy}

An efficient optimization strategy can decrease the complexity induced by the optimization algorithm. One possible approach is nested optimization that divides the problem into two optimization levels to exploit the linear-nonlinear structure of the problem (\cite{FAKIH20231450}). Furthermore, decomposition-based optimization methods are frequently employed to tackle large and complex optimization tasks by breaking a problem into smaller, more manageable subproblems. Well-known examples include (\cite{MITRAI2024108686}): Benders and Generalized Benders Decomposition, Lagrangian decomposition, the Alternating Direction Method of Multipliers (ADMM), and cross decomposition. Although these approaches have achieved broad success including in energy system optimization (\cite{LIU202443}), their advantage over solving the problem as a single “monolithic” problem is not always guaranteed. Moreover, their practical implementation, especially in real-time contexts, is not practical because of the numerous steps involved in their construction and their associated computational complexity.
\subsection{Optimization under uncertainties} \label{sec:optim_uncertain}
\subsubsection{Optimization under Aleatoric uncertainties}
The management of uncertainties has become an increasingly important topic in energy systems optimization. The reviews (\cite{fodstad_next_2022}) and  (\cite{lasemi_comprehensive_2022}) provide an overview of the various methods employed and their limitations. (\cite{lasemi_comprehensive_2022}) categorizes the methods into six groups, which can be ordered based on the level of available information about uncertainties (related to Walker's classification (\cite{walker_deep_2013})): Monte Carlo Simulation (MCS), Point Estimate Method (PEM), Scenario-Based Approaches (SBA) which include stochastic optimization (SO), Possibilistic/Fuzzy Approaches, Robust and Distributionally Robust Optimization (RO and DRO), Information Gap Decision Theory (IGDT), and Interval Approaches. Methods combining these techniques exists also.

According to the data presented by (\cite{lasemi_comprehensive_2022}), scenario-based approaches and RO methods are the most prevalent in the literature. (\cite{roald_power_2023}) focuses on these approaches for optimizing power systems, and present their characteristics in a more formal way. This is particularly supported by the formalism of Distributionally Robust Optimization or DRO (\cite{rahimian_distributionally_2022}), which integrates concepts of robust optimization, risk aversion, and chance-constrained optimization. The DRO approach helps mitigate the conservative tendencies of RO by allowing adjustments to the decision-makers' risk aversion levels and by allowing to model uncertainties even when no defined probability distribution is available.

Despite the mature theoretical framework and the important number of available methods, previous works highlight two major categories of unsolved cjallenges related to managing uncertainties. The first category of challenges pertains to the characterization of uncertainties (\cite{mavromatidis_review_2018}). While probability distributions can be constructed for some uncertainties, others are difficult or impossible to quantify, such as future energy prices or the likelihood of extreme events. Moreover, time-dependent uncertainties, for which generating sufficiently diverse and statistically coherent scenarios, remain challenging. The second category of challenges involves the computational time associated with these methods, which significantly increases compared to purely deterministic approaches. Many authors (\cite{roald_power_2023},  \cite{hoffmann_review_2024}, among others) highlight the difficulty of scaling and the trade-off between model detail and computation time. Unlike deterministic calculations with fixed parameters, methods involving uncertainties not only increase the number of necessary calculations but can also lead to infeasible cases or significantly longer resolution times for some parameter configurations.

%
Finally, data-driven methods for optimization under uncertainty exist and are summarized in (\cite{NING2019434}). These methods employ real data, rather than purely assumed models or distributions to model uncertain parameters. This makes solutions more robust to inaccuracies in probability assumptions. By leveraging machine learning techniques (e.g., clustering, density estimation, or deep generative models), these approaches learn richer representations of uncertainty and integrate them directly into the usual mathematical programming frameworks presented previously.

\subsubsection{Optimization under epistemic uncertainties}

Epistemic uncertainties can be reduced if additional data is collected and used to calibrate the model prior to the optimization process. \cite{jung_statistical_2022} provides a brief review of existing methods for parameter estimation under epistemic uncertainties. Existing methods include Kennedy and O’Hagan (KOH) framework (\cite{kennedy_bayesian_2001}), Optimization-based model calibration (OBMC) (\cite{lee_review_2019}). However, epistemic uncertainty cannot be completely eliminated in real world applications because of lack of data. In this case, since no probability functions are known for them, those are frequently modeled using intervals or, alternatively, using probability boxes (p-boxes) and the most suitable optimization frameworks are then non-probabilistic such as minimax or robust optimization. However, some contributions proposed to adapt probabilistic methods to this case, such as Bayesian optimization (\cite{celorrio_reliability-based_2021}). This implies mainly two important limitations: the difficulty of defining an appropriate uncertainty set and the prohibitive complexity when the size of uncertainty set or the number of uncertain parameters grows. 

\subsection{Prevalence of Data-Driven Techniques in Classical Approaches and their main limitations}
Many of the classical methods described in section \ref{sec:classic} inherently rely on data-driven techniques to process historical data or simulation snapshots for reducing complexity or managing uncertainty. In time-series reduction and spatial aggregation, clustering methods (e.g., k-means, hierarchical approaches)  compresses temporal information by identifying typical periods and reduce spatial complexity by aggregating regions with similar demand or resources profile. This speeds up subsequent optimization. However, these methods can be suboptimal if they fail to integrate physical system properties that are not fully captured by time-series data alone. In model Reduction, complex, high-fidelity models can be reduced to surrogate or reduced-order representations built from simulation or operational data. These data-driven substitutes replicate the main behavior of the system but require extensive offline sampling, and can lose accuracy or interpretability if applied to conditions beyond their training range. In scenario generation for optimization under incertainties, methods like generative machine learning use probabilistic parameters or empirical distributions to account for renewables variability, load uncertainty, or price volatility. While this captures some uncertainty, it relies heavily on the completeness and quality of available data and can become computationally burdensome as the number of scenarios grows.

Based on previous analysis, the high potential of data-driven models for reducing complexity and handling uncertainties is clearly demonstrated, these data-driven methods rel. However, these methods rely heavily on the quantity, the quality and completeness of available data. This can lead to either to overfitting or poor generalization. Moreover, purely empirical models tend to have limited physical interpretability, making decisions understanding and adoption very difficult. Finally, the computational burden of generating or handling large datasets is a serious limitation.

\section{Opportunity for hybrid modeling to mitigate these limitations}
To mitigate the limitations of data-driven methods, hybrid modeling approaches are emerging. Hybrid modeling mixes mechanistic (physics-based) knowledge with data-driven (ML) components. Hybrid approaches combine the strengths of both worlds: it maintains the interpretability and reliability offered by engineering and physical knowledge, while harnessing pattern-recognition capabilities and computational efficiency of data-driven methods. These techniques are applicable at different levels:

\textbf{1. Scenarios  generation: }
The generation of realistic input scenarios for optimization in energy systems can be achieved using hybrid generative models that incorporate both sparse data and physical knowledge. An illustrative example could be using a mechanistic model for global warming along with historical data to generate future weather conditions and demand curves. Some examples in the literature exist such as \cite{Afra_2023} who developed a physics-informed generative model for scenario generation of residential demand and \cite{kabir_2021} who proposed a hybrid model for accurate solar PV generation estimation. Hybrid models not only leverage the strengths of mechanistic understanding to impose physically realistic constraints, but also enhances the representation of uncertainty in the generated scenarios.

\textbf{2. Model reduction and calibration:} this aims to replacing the system's model by a hybrid model which reduces its complexity and facilitates the calibration process. Common strategies to develop a hybrid model include (\cite{rudolph_hybrid_2024,shah_hybrid_2025}): 
\textbf{- Series or Cascaded} One model feeds its outputs into the other (e.g., physics → ML or ML → physics). This was successfully applied to packed-bed heat storage modeling (\cite{PADMANABHAN2024113068}).
\textbf{- Parallel or Residual :} Physics and ML both run using the same inputs, and their results are merged.
\textbf{- Nested or Iterative Feedback} Physics and ML update one another over time steps, ideal for real-time control or dynamic processes. This approach was successfully applied to heat pumps modeling (\cite{mabrouk_neural-accelerated_2025}).
\textbf{- Physics-Informed Neural Networks (PINNs):} The neural network’s training loss enforces physical conservation laws, ensuring plausible physical behavior. This approach was succussfuly applied to AC power flows (\cite{JALVING2024109741}) or solar power plants (\cite{OSORIO2025119542}).

In addition to model reduction, hybrid models allow for more efficient calibration processes through the reduction of epistemic uncertainties. \cite{nair_e-pinns_2025} proposes a physics-informed machine learning strategy capable of quantifying epistemic uncertainty. Moreover, hybrid models allow easier integration in existing optimization under uncertainties frameworks. For example,  the posterior distributions of Bayesian hybrid models provide a principled
uncertainty set for stochastic programming, chance-constrained optimization, or robust optimization (\cite{eugene_learning_2023}).
Developed hybrid model can be wrapped in optimization problem solved using traditional solvers to provide final decisions.

\textbf{3. Optimization:} Hybrid methods can be used to directly predict optimal decision or to guide traditional solvers to accelerate their convergence. For instance, through an objective function that integrates governing laws, constraints, and goals, PINNs
enables top-down searches for optimal solutions (\cite{seo_solving_2024}). This approach was applied successfully to geothermal wells control (\cite{YAN2024123179}). Similarly, physics-informed reinforcement learning integrates domain-specific constraints within model-free policy learning to guide exploration and focus on feasible decision spaces. This was applied to both optimal power flow with distributed generation (\cite{Wu_2024}) and control of heating systems of buildings (\cite{saeed_physics-informed_2023}). 

\section{Conclusion}

This work underscores the increasing importance of advanced modeling and optimization methods as multi-energy systems grow in scale and technical sophistication. Conventional model-based optimization faces escalating computational requirements, especially when seeking high levels of temporal, spatial, and technological detail. Parallel to these challenges, the integration of uncertainties—both aleatoric and epistemic—complicates decision-making at every level, from real-time control to long-term planning. While classic strategies (e.g., time-series reduction, model order reduction, scenario generation) have proven effective, they often force trade-offs between accuracy and tractability.

Hybrid modeling offers a promising path forward by bridging the gap between purely mechanistic and purely data-driven techniques. Incorporating the interpretability and physical consistency of mechanistic models with the adaptability and speed of machine learning opens the door to improved accuracy, scalability, and uncertainty handling. Future developments in MES research will likely focus on refining these hybrid methods, including automating the process of selecting suitable model abstractions, developing robust frameworks for scenario generation under limited data, and better integrating uncertainty quantification into optimization. Ultimately, leveraging hybrid modeling approaches will be critical for designing and operating MES that can meet sustainability goals while maintaining reliability, cost-effectiveness, and adaptability to an evolving energy landscape.

\section*{REFERENCES}

\printbibliography[heading=none]

\section*{ACKNOWLEDGEMENT}
This work has been founded by the National French Research Agency (ANR) under grant number 22-PETA-0002. The HyMES project is part of the PEPR TASE research program.

\end{document}